\documentclass[runningheads,a4paper]{llncs}

\usepackage{graphicx,amssymb,mathtools, dsfont, amsfonts}
\usepackage{graphicx}
\usepackage{tikz}
\usepackage{tabularx}
\usepackage{pgfplots}

\usepackage{url}
\urldef{\mailsa}\path|{johannes.berger, frank.lenzen}@iwr.uni-heidelberg.de|
\urldef{\mailsb}\path|{neufeld, becker, schnoerr}@math.uni-heidelberg.de|

\newcommand{\keywords}[1]{\par\addvspace\baselineskip
\noindent\keywordname\enspace\ignorespaces#1}

\allowdisplaybreaks[4]

\DeclareTextFontCommand{\emph}{\em}

\begin{document}




\newcommand{\bitem}{\begin{itemize}}
\newcommand{\eitem}{\end{itemize}}
\newcommand{\mc}[1]{\mathcal{#1}}
\newcommand{\mb}[1]{\mathbb{#1}}
\newcommand{\mf}[1]{\mathfrak{#1}}
\newcommand{\ms}[1]{\mathscr{#1}}
\newcommand{\on}[1]{\operatorname{#1}}
\newcommand{\II}{\mathbb{I}}
\newcommand{\N}{\mathbb{N}}
\newcommand{\R}{\mathbb{R}}
\newcommand{\C}{\mathbb{C}}
\newcommand{\F}{\mathcal{F}}
\newcommand{\B}{\mathbb{B}}
\newcommand{\U}{\mathbb{U}}
\newcommand{\EE}{\mathbb{E}}
\newcommand{\V}{\mathbb{V}}
\newcommand{\Q}{\mathbb{Q}}
\newcommand{\Z}{\mathbb{Z}}
\newcommand{\PP}{\mathbb{P}}
\newcommand{\TT}{\mathbb{T}}
\newcommand{\bpm}{\begin{pmatrix}}
\newcommand{\epm}{\end{pmatrix}}
\newcommand{\bsm}{\left(\begin{smallmatrix}}
\newcommand{\esm}{\end{smallmatrix}\right)}
\newcommand{\T}{\top}
\newcommand{\ul}[1]{\underline{#1}}
\newcommand{\ol}[1]{\overline{#1}}
\newcommand{\la}{\langle}
\newcommand{\ra}{\rangle}
\newcommand{\si}{\sigma}
\newcommand{\SI}{\Sigma}
\newcommand{\mrm}[1]{\mathrm{#1}}
\newcommand{\msf}[1]{\mathsf{#1}}
\newcommand{\mfk}[1]{\mathfrak{#1}}
\newcommand{\row}[2]{{#1}_{#2,\bullet}}
\newcommand{\col}[2]{{#1}_{\bullet,#2}}
\newcommand{\df}[2]{\frac{\partial #1}{\partial #2}}
\newcommand{\p}{\partial}
\newcommand{\veps}{\varepsilon}
\newcommand{\toset}{\rightrightarrows}
\newcommand{\w}{\omega}
\newcommand{\gdw}{\Leftrightarrow}
\newcommand{\vphi}{\varphi}
\newcommand{\ora}[1]{\overrightarrow{#1}}
\newcommand{\ola}[1]{\overleftarrow{#1}}
\newcommand{\oset}[2]{\overset{#1}{#2}}
\newcommand{\uset}[2]{\underset{#1}{#2}}
\newcommand{\SE}{\operatorname{SE}_{3}}
\newcommand{\se}{\mathfrak{se}_{3}}
\newcommand{\SO}{\operatorname{SO}_{3}}
\newcommand{\so}{\mathfrak{so}_{3}}
\newcommand{\Ad}{\operatorname{Ad}}
\newcommand{\dist}{\operatorname{dist}}
\newcommand{\etr}{\operatorname{etr}}
\newcommand{\vex}{\operatorname{vex}}
\newcommand{\Psym}{\mathbb{P}_{s}}
\newcommand{\Pskew}{\mathbb{P}_{a}}
\newcommand{\vecso}{\operatorname{vec}_{\mathfrak{so}}}
\newcommand{\vecse}{\operatorname{vec}_{\mathfrak{se}}}
\newcommand{\matso}{\operatorname{mat}_{\mathfrak{so}}}
\newcommand{\matse}{\operatorname{mat}_{\mathfrak{se}}}
\newcommand{\Hess}{\operatorname{Hess}}
\newcommand{\grad}{\operatorname{grad}}
\newcommand{\tr}{\operatorname{tr}}
\newcommand{\argmin}{\operatorname{arg\,min}}
\newcommand{\kronse}{\otimes_{\mathfrak{se}}}
\newcommand{\diag}{\operatorname{diag}}
\newcommand{\Exp}{\operatorname{Exp}}
\newcommand{\D}{\mathbf{D}}

\newcommand{\eins}{\mathds{1}}

\newcommand{\cperp}[3]{{#1} \perp\negthickspace\negthinspace\negthickspace\perp {#2} \,|\, {#3}}

\newcommand{\LG}[1]{\mathrm{#1}}
\newcommand{\Lg}[1]{\mathfrak{#1}}

\newcommand{\st}[1]{{\scriptstyle #1}}
\newcommand{\sst}[1]{{\scriptscriptstyle #1}}


\mainmatter  %

\title{Second Order Minimum Energy Filtering on $\SE$ \\ with Nonlinear Measurement Equations}

\titlerunning{Minimum Energy Filtering on $\SE$ \\ with Nonlinear Measurement Equations}
\authorrunning{J.~Berger \and A.~Neufeld \and F.~Becker \and F.~Lenzen \and C.~Schn\"orr}
\author{Johannes Berger \and Andreas Neufeld \and Florian Becker \and \\ Frank Lenzen \and Christoph Schn\"orr}
\institute{IPA \& HCI, University of Heidelberg, Germany \\ \mailsa\\ \mailsb}

\maketitle

\begin{abstract}
Accurate camera motion estimation is a fundamental building block for many Computer Vision algorithms.
For improved robustness, temporal consistency of translational and rotational camera velocity is often assumed by propagating motion information forward using stochastic filters.
Classical stochastic filters, however, use linear approximations for the non-linear observer model and for the non-linear structure of the underlying Lie Group~$\SE$ and have to approximate the unknown posteriori distribution. 
In this paper we employ a non-linear measurement model for the camera motion estimation problem that incorporates multiple observation equations. We solve the underlying filtering problem using a novel Minimum Energy Filter on $\SE$ and give explicit expressions for the optimal state variables. Experiments on the challenging KITTI benchmark show that, although a simple motion model is only employed, our approach improves rotational velocity estimation and otherwise is on par with the state-of-the-art.

\end{abstract}

\section{Introduction}\label{sec1}

Camera motion estimation
is an important task in autonomous driving for which the ego-motion of the camera is fully determined by images from cameras mounted on the car. 
Most approaches require only temporal correpondences~\cite{Becker2013} or additional depth information~\cite{geiger2011stereoscan,kitt2010visual} e.g.\ obtained from stereo estimation.
Given two frames 
and a depth map, the underlying motion of the camera can be determined uniquely. 
However, two-frame methods are sensitive to noise and thus past information needs to be propagated with filtering approaches. 
Stochastic filters require assumptions about the a posteriori distribution, which is often unknown and thus can not be modeled adequately. Furthermore, an adaption to Lie groups is unknown for almost all stochastic filters.
Application of state-of-the-art particle filters is limited due to the high amount of required particles.
Mortensen~\cite{Mortensen1968} derived a second order deterministic filter for the classical filtering problem on~$\R^{n}$ based on classical control theory. This result has been generalized to Lie groups~\cite{saccon2013second}.

In this article we present a filtering model with state and observation equations for the camera motion estimation problem. We adapt the approach~\cite{saccon2013second} to this model on~$\SE$ and generalize it to incorporate \emph{multiple} measurement equations depending non-linearly on the camera motion. We also show how the abstract exponential functor~\cite{saccon2013second} can be computed explicitly and  finally derive a matrix representation of the inverse Hessian operator on~$\se$ using special Kronecker products. Numerical experiments show the fast convergence of the filter. 

\subsubsection{Related work.} The task of estimating the current state of a dynamical system only based on past observations of the system 
is known as \emph{filtering}. 
In the last century numerous \emph{stochastic filters} have been developed starting from the seminal work of Kalman~\cite{kalman1960new}. See \cite{Stochastic-Filtering-09} for an overview and background. Since for non-linear dynamical systems with non-Gaussian noise processes this problem cannot be solved exactly, several approaches tried to cope with these non-linearities~\cite{Daum2005}. In general, as the a-posteriori distribution is unknown, most Gaussian filters are doomed to fail. State-of-the-art particle filters~\cite{Daum2003Curse} alleviate the problem of not knowing the distribution. For large dimensions, however, generating enough particles becomes infeasible. Brigo et al.~\cite{NonlFilterExpFamily-99} use exponential families to model the a posteriori distribution, but the choice of an admissible sufficient statistic is critical and too restricted for our multiple measurement model.

For the considered filtering problems on Lie groups we have to take  into account the non-linear geometry of the manifold to find an optimal filter. Markley~\cite{markley2003attitude} worked out a method for filtering problems on~$\SO$ whereas~\cite{kwon2007} investigates particle filters on~$\SE$. While the dimension of the embedding space is not excessively large (e.g, 16 for $\SE$), the generation of samples \emph{on} the respective Lie group is considerably more expensive.
Mortensen~\cite{Mortensen1968} derived a \emph{deterministic and recursive second order optimal filter} based on results of \emph{control theory} and the \emph{dynamic principle}. In the last years this approach was generalized to specific Lie groups~\cite{Zamani2012,saccon2013second}. In various scenarios it has been shown that minimum energy filters have an \emph{exponential} convergence rate~\cite{krener2003convergence} and perform superior to extended Kalman filters~\cite{Zamani2012}. 

\subsubsection{Contributions.}
\begin{itemize}
\item Formulation of \emph{filtering equations} for the camera motion estimation problem on rigid scenes with constant motion assumption; 
\item adaptation of the second order Minimum Energy Filter~\cite{saccon2013second} such that it incorporates \emph{multiple} and \emph{non-linear} measurement equations;
\item derivation of \emph{explicit} ordinary differential equations of the optimal state and the inverse Hessian for which we derive a matrix representation;
\item experiments that show the comparable performance in accuracy of the camera motion against a state-of the art-method~\cite{geiger2011stereoscan} on the challenging real-life KITTI benchmark.
\end{itemize}

\subsubsection{Preliminaries.}
We use the notation~$\SE$ to denote the special Euclidean group equipped with its tangent space~$T_{E}\SE$  and Lie algebra~$\se$. 
Tangent vectors $E\Gamma \in T_{E}\SE$ are obtained as evaluations of left-invariant vector fields $\D L_{E}(I)[\Gamma]$ that one-to-one correspond to the  tangent vector $\Gamma \in \se$. Here, $L_{E}F = E F$ denotes left-translation and $\D L_{E}$ denotes the differential of $L_{E}$.
$\SE$ can be identified with a matrix Lie subgroup of $\on{GL}(4)$. We adopt the Riemannian metric 
$\la X, Y \ra_{E} := \la E^{-1} X, E^{-1} Y \ra, \forall X, Y \in T_{E} SE_{3}$,
where $\la A, B \ra = \tr(A^{\T} B)$ denotes the canonical matrix inner product.
The Riemannian gradient~$\grad f$ of a differentiable function $f:\SE \rightarrow \R$ is defined through the directional derivative $\mathbf D f(E)[E\Omega] =: \la \grad f(E), E\Omega  \ra_{E}$ for all $\Omega \in \se$. The Riemannian Hessian~$\Hess f (E)[\cdot]: T_{E}\SE \rightarrow T_{E}\SE$ at~$E \in \SE$ is defined through the relation~$\la \Hess f(g)[E\Gamma], E\Omega \ra := \D( \D f(E) [E\Omega])[E\Gamma] - \D f(E)[\nabla_{E\Gamma} E\Omega]$ for all $\Gamma,\Omega\in \se$. Here, $\nabla$ denotes the Riemannian (Levi-Civita) connection.
The Lie algebra~$\se$ can be associated with a 6-dimensional vector space and we define the operation~$\vecse: \se \rightarrow \R^{6}$ given by
\begin{equation}
\vecse \left(\left(\begin{smallmatrix} 
0 & -\gamma_{3} & \gamma_{2} & \gamma_{4} \\
\gamma_{3} & 0 & -\gamma_{1} & \gamma_{5}\\
-\gamma_{2} & \gamma_{1} & 0 & \gamma_{6} \\
0 & 0 & 0 & 0 
 \end{smallmatrix}\right)\right) = (\gamma_{1},\gamma_{2}, \gamma_{3}, \gamma_{4},\gamma_{5}, \gamma_{6})^{\T}.
\end{equation}
We denote the inverse operation by $\matse: \R^{6} \rightarrow \se$ and the projection onto the Lie algebra by $\on{Pr}:\on{GL}(4) \rightarrow \se$. The standard basis of $\se$ is given by $\{B_{i} = \matse(b_{i})\}$, where $b_{i}, i=1,\dots,6$ is the standard basis of $\R^{6}.$

\section{Minimum Energy Filtering Approach} The classical filtering problem consists of a \emph{state equation} that describes the dynamics of an \emph{unknown} state $E(t)$ and \emph{observation equations} connecting measurements to the state of the system. These \emph{real-valued} equations are given by
\begin{align}
\dot{E}(t) & = f_{t}(E(t)) + \delta(t),\quad E(0) = E_{0,}& \mbox{(state)} \label{eq:general-state} \\
y(t) & = h_{t}(E(t)) + \epsilon(t), & \mbox{(observation)}\label{eq:general-observation}
\end{align}
where the functions $f_{t}$ and $h_{t}$ describe the state and observation dynamics, respectively, and $\delta(t), \epsilon(t)$ are noise processes. Stochastic filters usually understand these equations as stochastic differential equations and try to find for each~$t$ the maximum of the a posteriori distribution~$P(E(t) \vert y(s), s \leq t)$. In contrast, Mortensen~\cite{Mortensen1968} investigated~\eqref{eq:general-state} and~\eqref{eq:general-observation} from control theory point of view: He considered the equation~\eqref{eq:general-state} as a dynamical system, controlled by a control process~$\delta(t)$ such that the residual~$\lVert \epsilon(t) \rVert = \lVert y(t) - h_{t}(E(t))\rVert$ is minimized.

\subsection{Measurement Model for Ego-Motions}
\subsubsection{Supported Models.}
In this section we derive an optical flow observer model with corresponding state equation on the Lie group $\SE.$ They support two models that incorporate two different kinds of given data, i.e.
\begin{itemize}
\item temporal optical flow and stereo matches (stereo approach),
\item temporal optical flow and depth map (monocular approach).
\end{itemize}
Since the proposed filter supports both models and shares the same derivation, we only consider the \emph{monocular} model in the following.

\subsubsection{Optical Flow Induced by Egomotion.}
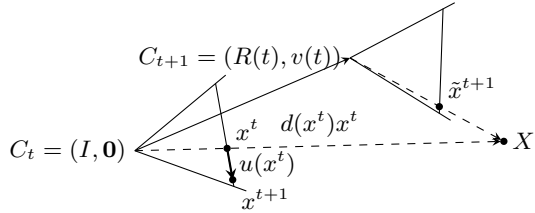
\begin{figure}[htbp]
\begin{center}
\begin{center}
\begin{tikzpicture}[scale=1.2,rotate=-80]
\def\camera(#1,#2,#3){
\draw[thin, shift={(#1,#2)},rotate=30+#3]   (0,0) -- ++(0,1.3);
\draw[thin, shift={(#1,#2)},rotate=-30+#3]  (0,0) -- ++(0,1.3);
\draw[thin, shift={(#1,#2)},rotate=#3]      (-0.58,1.) -- ++(1.16,0.);
}
\def\r{0.04};

\def \hone{-0.6};
\def \htwo{2.5};

\coordinate (origin) at (0,0);
\coordinate (Cr) at (2.0, 0.);
\coordinate (X) at (0.6,4);
\coordinate (h) at (\hone,\htwo);
\coordinate(x0l) at (0.15,1);
\coordinate(x1l) at (0.502069, 1.); 
\coordinate(x1ltilde) at (0.1,3.37);

\draw node [anchor=west] at (X) {$X$};
\draw node [anchor=south west] at (x0l) {$x^{t}$};
\draw node [anchor=east] at (origin) {$C_{t} = (I,\mathbf{0})$};
\draw node [anchor=south, below, left] at (h) {$C_{t+1} = (R(t),v(t))$};
\draw node [anchor=north west] at (x1l) {$x^{t+1}$};
\draw node [anchor=south west] at (x1ltilde) {$\tilde{x}^{t+1}$};

\camera(0,0,0);
\camera(\hone,\htwo,-12);

\fill (x0l) circle (\r);
\fill (x1l) circle (\r);
\fill (x1ltilde) circle (\r);

\fill (X) circle(\r);
\draw[shorten >=0.04cm,->,thin,>=stealth, dashed] (0,0) -- (X) node[midway, above] {$d(x^{t}) x^{t}$};
\draw[shorten >=0.04cm,->,thin,>=stealth, dashed] (h) -- (X);
\draw[->,>=stealth,line cap=triangle 90] (0,0) -- (h) node[midway, left] {};
\draw[->, thick,>=stealth,line cap=triangle 90] (x0l) -- (x1l) node[midway, right] {$u(x^{t})$};
\end{tikzpicture}
\end{center} 
\caption{Setup for temporal optical flow for either a given depth map or  stereo matchings. Correspondences are given by $x^{t+1}= x^{t} + u(x^{t,l}),$ and $\tilde{x}^{t+1} = v(t) + R(t)x^{t+1}$ denotes the perspective projection of~$X$ to the camera plane of $C_{t+1}.$} 
\label{fig:camera-model}
\end{center}
\end{figure}

In this work we denote the time space by~$T := \R_{\geq0}$ and the image sequence recorded by a camera moving through a static scene by~$f = \{f_{t}, t\in T\}$.
Let~$(R(t),v(t)) = E(t) \in \SE$ be the incremental camera motion from frame $f_{t}$ to $f_{t+1}$. 
W.l.o.g.\ at time~$t$ we set the coordinate system to be identical to the one of the camera recording~$f_t$, i.e.\
the extrinsic camera parameters are~$C_{t} = (I, 0)$ and $C_{t+1} = (R(t),v(t)) $. 
Let~$X \in \R^{3}$ be a scene point and we denote its perspective projection into camera~$C_t$ by
$x^{t}= (x_{1}^{t},x_{2}^{t},1)^{\T}=P_C(X)$ where $P_{C}((x_{1},x_{2},x_{3})^{\T}) := x_{3}^{-1}(x_{1},x_{2},x_{3})^{\T}$. Furthermore, we denote by~$d(x^{t}) \in \R$ the depth of~$X$ in camera~$C_t$ such that
the scene point can be reconstructed by~$X = x^{t} d(x^{t})$,
see Fig.~\ref{fig:camera-model}. 

With this relation we obtain the optical flow induced by the camera motion~$(R(t),v(t))$ and depth~$d(x^t)$: 
\begin{align}
u(R(t),v(t),x^{t},d(x^{t})) + x^{t}=  & P_{C}\Bigr(R(t)^{\T}( x^{t} d(x^{t})-v(t))\Bigr)\label{eq:flow_observer_depth}.
\end{align}
We introduce index~$k$ to distinguish multiple observations and reformulate the observer equations in terms of~$E(t)=(R(t),v(t)) \in \SE$ yielding
\begin{equation}\label{eq:observation-dynamics-flow}
u(E(t),g_k(t)) + x_k^{t} = P_{C}(\hat{I}E(t)^{-1}g_k(t)), \quad\hat{I} := \left(\begin{smallmatrix}1 & 0 & 0 & 0 \\ 0 & 1 & 0 & 0 \\ 0 & 0 & 1 & 0 \end{smallmatrix}\right)
\end{equation}
where~$g_k(t) = g(t;x_k^{t}, d(x_k^{t})) :=  ( d(x^{t}_k) ({x_k^{t}})^{\T}  ,1)^\T \in \R^{4}$
for a single measurement.

\section{Minimum Energy Filter Derivation}
\subsubsection{State and observation equations.}
As we want to recover the camera motion from the image data, we cannot incorporate any prior knowledge of the cameras' kinematics such as data from external acceleration sensors. The only assumption we make is to demand a constant camera motion $E(t) \in \SE$ that is influenced by a noise process $\delta(t) \in \se$, which also models accelerations. This kinematic state equation on $\SE$ without dynamics $f_{t}$ (i.e. $f_{t}\equiv 0$) 
is given by
\begin{equation} \label{eq:state}
\dot{E}(t) = E(t) \delta(t) \in T_{E(t)}\SE, \quad E(0) = E_{0} \in \SE.
\end{equation}
We incorporate multiple flow observations $y_{k}(t):=u_{k}(E(t), g_{k}(t),t) + x_{k}$ at different image points $x_{k}$ for $k\in \{1,\dots, n\}$ that depend on the ego-motion $E(t)$ and are corrupted by noise vectors $\epsilon_{k}(t) \in \R^{3}.$ This gives
\begin{equation}\label{eq:observer}
y_{k}(t) = h_{k}(E(t)) + \epsilon_{k}(t)\,, \quad  k\in \{1,\dots,n\}\, .
\end{equation}
Here we used the observation functions $h_{k}:\SE \times \R^{4} \rightarrow \R^{3}, k \in \{1,\dots,n\}$ that we define as is Eq.~\eqref{eq:observation-dynamics-flow}, i.e. $h_{k}(E(t)) = h_{k}(E(t),g_{k}(t)):= P_{C}(\hat{I}E(t)^{-1} g_{k}(t))$.

\subsubsection{Energy function.} We want to find the camera motion that describes the observation process best up to a small error $\epsilon$, i.e. we want to minimize the residual  $\lVert \epsilon_{k}(t) \rVert_{Q} ^{2}$ ($\lVert \cdot\rVert_{Q}^{2}:= \cdot^{\T}Q\cdot$) for all $t$ such that the dynamical system on $E(t)$ \eqref{eq:state} is also fulfilled. The latter means that also the error term $\lVert \vecse(\delta(t))\rVert_{S}^{2}$ is minimized. Here $S \in \R^{6\times 6}$ and $Q \in \R^{3\times 3}$ are symmetric and positive definite weighting matrices. We define the following energy function
\begin{equation}\label{eq:cost-function}
\mathcal{J}(\epsilon, \delta, t_{0}, t):= m_{0}(E_{0}, t,t_{0}) + \int_{t_{0}}^{t}  c(\delta, \epsilon,\tau,t)  \, d\tau,  
\end{equation}
where $\epsilon := \{\epsilon_{k}, k=1,\dots,n \}$, $\delta := (\delta(\tau),\tau \leq t), t \in T$, $c: \se \times \R^{3n}  \times T\times T \rightarrow \R$ is a quadratic penalty function for $\delta$ and $\epsilon$ given by
\begin{equation}\label{eq:penalty-function-c}
c(\delta, \epsilon, \tau,t) := \tfrac{1}{2}e^{-\alpha(t-\tau)}\Bigl( \lVert \vecse(\delta(\tau))\rVert_{S}^{2} + \sum_{k=1}^{n} \lVert \epsilon_{k}(\tau)\rVert_{Q}^{2} \Bigr)
\end{equation}
and $m_{0}: \SE \times T \times  T \rightarrow \R_{\geq 0}$, $(E_{0}, t_{0},t_{1}) \mapsto  \tfrac{1}{2} e^{-\alpha(t-t_{0})}\tr((E_{0}-\eins)^{\T}(E_{0}-\eins))$, with $\eins$ being the identity matrix in $\R^{4\times 4}$,
is a penalty function for the initial condition.
Here we also used the idea of a decay rate $\alpha \geq 0$ from \cite{saccon2013second} at which old information is forgotten. 
To incorporate the observations \eqref{eq:observer} we substitute the error term $\epsilon_{k}(t)= \epsilon_{k}(E(t),t)$ by 
$y_{k}(t) - h_{k}(E(t))$ in Eq.~\eqref{eq:penalty-function-c}.

\subsubsection{Optimal control problem.}

The optimal control theory allows us to determine the optimal control input $\delta$ that minimize the energy $\mathcal J(\epsilon(E(t),t),\delta, t_{0},t)$ for each $t \in T$ subject to the state constraint \eqref{eq:state}. To be precise, we want to find $\delta\vert_{[t_{0},t]}$ for all $t\in T$ and fixed $E \in \SE$ defining the value function
\begin{equation}\label{eq:val_fun}
\mathcal V(E(t),t) := \min_{\delta \vert _{[t_{0},t]}} \mathcal J(\epsilon(E(t),t),\delta, t_{0},t) \,\, s.t.\,\dot{E}(t) = E(t) \delta(t),\, E(0)=E_{0}\, .
\end{equation}
The optimal trajectory is ${E^{\ast}}(t) := \argmin_{E(t) \in \SE} \mathcal{V}(E(t),t)$ for all $t \in T$ and $\mathcal{V}(E,t_{0}) = m_{0}(E_{0},t_{0},t_{0}).$ This problem is a classical optimal control problem, for which the classical Hamilton-Jacobi theory \cite{jurdjevic1997,athans1966optimal} gives the well known Hamilton-Jacobi-Bellman equation. Pontryagin \cite{athans1966optimal} proved that the minimization of the Hamiltonian gives a solution of the corresponding optimal control problem (\emph{Pontryagin's Minimum Principle}).

However, since $\SE$ is a non-compact Riemannian manifold we cannot apply the classical Hamilton-Jacobi theory for real-valued problems (cf.~\cite{athans1966optimal}).
Instead we follow the approach of Saccon et al.~\cite{saccon2013second} which derive a \emph{left-trivialized optimal Hamiltonian} based on control theory on Lie groups \cite{jurdjevic1997}, which is given by $\tilde{\mathcal H}^{-}: \SE \times \se \times \se \times T \rightarrow \R,$
\begin{equation}\label{eq:preHamiltonian}
\tilde{\mathcal{H}}^{-}(E,\mu, \delta, t) = c(\delta, \epsilon(E,t), t_{0},t) - \la \mu, \delta \ra\,.
\end{equation}
The minimization of \eqref{eq:preHamiltonian} w.r.t. the variable $\delta$ leads \cite[Prop.~4.2]{saccon2013second} to the optimal Hamiltonian $\mathcal{H}^{-}(E,\mu,t):= \tilde{\mathcal{H}}^{-}(E,\mu,\delta^{\ast},t),$ where $\vecse(\delta^{\ast})= e^{\alpha(t-t_{0})} S^{-1} \vecse(\mu).$ $\mathcal H^{-}$ is given by
\begin{equation}\label{eq:left-trivial-hamilton}
\mathcal{H}^{-}(E, \mu, t) = \tfrac{1}{2}e^{-\alpha(t-t_{0})} \sum_{k=1}^{n} \lVert y_{k}-h_{k}(E)\rVert_{Q}^{2}- \tfrac{1}{2}e^{\alpha(t-t_{0})} \la  \mu, \matse(S^{-1}\vecse(\mu)) \ra  \,.
\end{equation}
In the next section we compute explicit ordinary differential equations for the optimal state ${E^{\ast}}(t)$ for each $t\in T$ that consists of different derivatives of the left trivialized Hamilton function \eqref{eq:left-trivial-hamilton}.

\subsection{Recursive Filtering Principle of Mortensen}

In order to find a recursive filter we compute the total time derivative of the optimality condition on the value function, which is
\begin{equation}
\grad_{1} \mathcal V({E^{\ast}},t) = 0,
\end{equation}
 for each $t \in T.$ This equation must be fulfilled by an optimal solution of the filtering problem ${E^{\ast}} \in \SE.$  Unfortunately, because the filtering problem is in general infinite dimensional, this leads to an expression containing derivatives of every order. In practice (cf. \cite{Zamani2012,saccon2013second}) derivatives of third order and higher are neglected, since they are complicated to compute. Omitting these leads to a second order approximation of the optimal filter. 
 
 \begin{theorem}
 \label{prop:main-result}
 The differential equations of the second order Minimum Energy Filter for our state \eqref{eq:state} and nonlinear observer \eqref{eq:observer} model is given by
 \begin{align}
\dot{{E^{\ast}}} = & - {E^{\ast}}\matse\bigl(P\vecse(\sum_{k}\on{Pr} \bigl(A_{k}({E^{\ast}}))  \bigr)), \quad E^{\ast}(t_{0}) = E_{0} \,, \label{eq:final-estimate-E}\\
\begin{split}
\dot{P} = & - \alpha P + S^{-1}-  P \sum_{k}\bigl(\tilde{\Gamma}_{\vecse(\on{Pr}(A_{k}({E^{\ast}})))}  + D_{k}({E^{\ast}}) \bigr) P  \\
& \quad - \tilde{\Gamma}_{\vecse({E^{\ast}}^{-1} \dot{{E^{\ast}}})}^{\ast}P + P (\tilde{\Gamma}_{\vecse({E^{\ast}}^{-1} \dot{{E^{\ast}}})}^{\ast})^{\T},  \quad P(t_{0}) = P_{0}\,,
\end{split} \label{eq:final-estimate-Q} 
\end{align}
where 
$A_{k}(E) = A_{k}(E,g_{k}):= 
\bigl(\kappa_{k}^{-1}\hat{I}  -  \kappa_{k}^{-2}  \hat{I}   E^{-1}  e_{3} g_{k}^\T   \hat{I} \bigr)^{\T}
 Q (y_{k} - h_{k}(E))  g_{k}^{\T}  E^{-\T}$,  $\kappa_{k}:=\kappa_{k}(E):= e_{3}^{\T} \hat{I}E^{-1}g_{k}$
and
$\on{Pr}(A):= \argmin_{\Omega \in \se} \la \Omega, A \ra$.
$D_{k}(E)$ is derived in the appendix in Eq.~\eqref{eq:def_Dk} and the matrix valued functions $\tilde{\Gamma}_{\cdot},\tilde{\Gamma}_{\cdot}^{\ast}: \R^{6} \rightarrow  \R^{6\times 6}$ come from the vectorization of the connection functions. Their components are given by $(\tilde \Gamma_{z})_{ij} := \sum_{k=1}^{6}\hat\Gamma_{jk}^{i} z^{k}$ and $(\tilde \Gamma_{z}^{\ast})_{ik} := \sum_{j=1}^{6}\hat\Gamma_{jk}^{i} z^{j}$ with $z \in \R^{6}$ and  Christoffel-Symbols $\hat{\Gamma}_{jk}^{i}:=\Gamma_{kj}^{i}$ from \cite{Zefran1999}.

The initial $P_{0}\in \R^{6\times 6}$ is given by $P_{0}\vecse(\Omega) =  \vecse\bigl((E \Hess m_{0}(E_{0})[E\Omega])^{-1}\bigr),$
 and $E_{0}$ is an initialization in $\SE.$
 \end{theorem}
The sketch of the proof is given at the end of this section. For the proof we need some lemmas listed below, the proofs of which can be found in the appendix.
 
We adapt the minimum energy filter for general Lie groups derived in \cite{saccon2013second} to our nonlinear measurement model on $\SE$. Following \cite[Eq. (37)]{saccon2013second} the estimate of the optimal state ${E^{\ast}}$ is given by
\begin{equation}\label{eq:ML-E}
\dot{{E^{\ast}}} = - {E^{\ast}}   Z({E^{\ast}},t)^{-1}\bigl(\grad_{1} \mathcal H^{-}({E^{\ast}},\mu,t)) \bigr),
\end{equation}
which contains the second order information matrix $Z({E^{\ast}},t):\se \rightarrow \se$ of the value function $\mathcal V$ (cf.~\eqref{eq:val_fun}), defined by
$Z({E^{\ast}},t)(\Omega) \mapsto {E^{\ast}}^{-1} \Hess_{1}\mathcal{V}({E^{\ast}},t)[{E^{\ast}}\Omega].$ The gradient of the Hamiltonian in \eqref{eq:ML-E} is given in the following lemma.
 \begin{lemma}\label{lemma:riemannian-gradient}
The Riemannian gradient on $T_{E}\SE$ of the Hamiltonian $\mathcal H^{-}(E, \mu,t)$ with $A_{k}(E) = A_{k}(E,g_{k})$ defined in Theorem \ref{prop:main-result} is given by
\begin{equation}\label{eq:explicit-grad-hamiltonian}
\grad_{1} \mathcal H^{-}(E, \mu, t) = e^{-\alpha(t-t_{0})} \sum_{k} E \on{Pr}(A_{k}(E)).\,
\end{equation}
 \end{lemma}
In order to derive a second order filter, we also need a recursive expression of the operator $Z$, that has been derived in  \cite[Eq. (51)]{saccon2013second} and  is approximately
\begin{align}
\begin{split} \label{eq:dynamic_programming_Hessian}
\hspace{-0.0cm}&\frac{d}{dt}  Z({E^{\ast}}(t),t) \approx  \omega_{{E^{\ast}}^{-1} \dot{{E^{\ast}}}}^{\ast} \circ Z({E^{\ast}},t) + Z({E^{\ast}},t) \circ \omega_{{E^{\ast}}^{-1} \dot{{E^{\ast}}}} \\
&\hspace{-0.2cm}  + Z({E^{\ast}},t) \circ \Hess_{2} \mathcal H^{-} ({E^{\ast}},0,t) \circ Z({E^{\ast}},t) + {E^{\ast}}^{-1} \Hess_{1} \mathcal H^{-}({E^{\ast}},0,t) \circ {E^{\ast}} ,\hspace{-0.0cm}
\end{split}
\end{align}
where third order derivatives are neglected. For the computation of the Hessian we need implicitly the Riemannian connection $\nabla$ with connection function $\omega_{\Omega}\Delta:=E\nabla_{\Omega}\Delta.$ The dual operator $\omega_{\Omega}^{\ast} \cdot$ is given by$\la \omega_{\Omega}^{\ast} \Delta, \Theta \ra := \la \Delta, \omega_{\Omega}\Theta \ra$.
 
 Next, we derive a matrix representation for all terms in Eq. \eqref{eq:dynamic_programming_Hessian} provided  by the $\vecse-$operation defined in section \ref{sec1} and the following lemmas.

 \begin{lemma}[Matrix representation of $Z$] \label{lem:2} Let $Z({E^{\ast}},t):\se \rightarrow \se$ be the operator in equation \eqref{eq:ML-E}. Then there exists a matrix $K = K(t)\in \R^{6\times 6}$ such that we can vectorize $Z({E^{\ast}},t)(\Omega)$ for each $\Omega \in \se$. Then it holds $\vecse(Z({E^{\ast}},t)(\Omega)) = K(t) \vecse(\Omega)$, and thus $\vecse(d/dt Z({E^{\ast}},t)(\Omega)) = \dot{K}(t)\vecse(\Omega),$ as well as
\begin{enumerate}
	\item $\vecse(\omega_{{E^{\ast}}^{-1}\dot{{E^{\ast}}}}^{\ast}  Z({E^{\ast}},t) \circ \Omega)  = (\tilde{\Gamma}_{\vecse({E^{\ast}}^{-1}\dot{{E^{\ast}}})}^{\ast})^{\T} K(t) \vecse(\Omega)$ \, \label{lem2-5}
	\item $\vecse(Z({E^{\ast}},t) \circ \omega_{{E^{\ast}}^{-1}\dot{{E^{\ast}}}} \Omega)  = K (t)\tilde{\Gamma}_{\vecse({E^{\ast}}^{-1}\dot{{E^{\ast}}})}^{\ast} \vecse(\Omega)$\, , \label{lem2-4}
		\item $\vecse( Z({E^{\ast}},t)  (\Hess_{2}  \mathcal H^{-}({E^{\ast}},0,t)[ Z({E^{\ast}},t) (\Omega)]) ) 
\\ \quad \quad=  - e^{\alpha(t-t_{0})} K(t)S^{-1}K(t) \vecse(\Omega)\, ,$ \label{lem2-3}
 \end{enumerate}
with $\tilde \Gamma_{\cdot}$ and $\tilde \Gamma_{\cdot}^{\ast}$ from Thm.~\ref{prop:main-result}. \end{lemma}
Finally we have to apply the $\vecse-$operation to the last remaining term in \eqref{eq:dynamic_programming_Hessian}:
\begin{lemma}\label{lem3} It holds
\begin{align*}
\vecse({E^{\ast}}^{-1}& \Hess_{1} \mathcal H^{-}({E^{\ast}},0,t) [{E^{\ast}}\Omega])  \\& \qquad = e^{-\alpha(t-t_{0})}   \sum_{k} \bigl( \tilde{\Gamma}_{\vecse(\on{Pr}(A_{k}({E^{\ast}})))}  + D_{k}({E^{\ast}}) \bigr)  \vecse(\Omega)
\end{align*}
where $D_{k}(\cdot):\SE \rightarrow \R^{6\times6}$ and $\tilde{\Gamma}_{\cdot}:\R^{6} \rightarrow \R^{6\times 6}$ are given in the appendix.
\end{lemma} 
\begin{figure}[htbp]
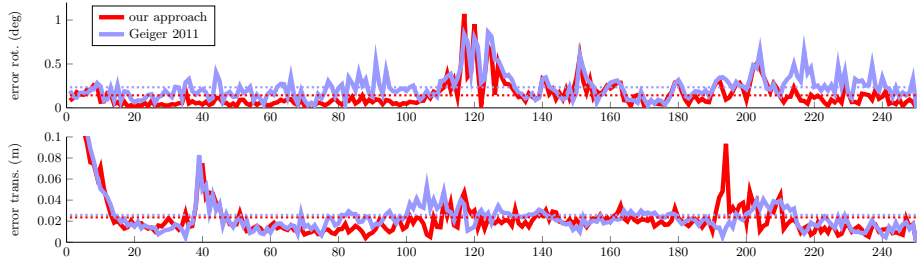

\begin{center}
 \input{errorRot}\\
 \input{errorTrans}\\
 \vspace{-1em}
\caption{Comparison of the rotational error (top, in degree) and the translational (bottom, in meters) of our approach and Geiger et al.~\cite{geiger2011stereoscan} on the first 250 frames of sequence 0 of the Kitti odometry benchmark. We used the parameters ($\alpha=2,\, S = \diag((s_{1},s_{1},s_{1},s_{2},s_{2},s_{2})), s_{1} = 10^{-3}, s_{2}=10^{-6}, Q = 0.02 \eins$). The dotted lines indicate the mean errors. 
In the translational part (bottom) both methods are competitive. In the rotational part (top) we outperform \cite{geiger2011stereoscan}.}
\label{fig:comparison-geiger}
\end{center}
\end{figure}
\renewcommand*{\proofname}{Sketch of proof of Thm.~\ref{prop:main-result}}
\begin{proof}
Insertion of \eqref{eq:explicit-grad-hamiltonian} into \eqref{eq:ML-E} and  $\vecse(Z({E^{\ast}},t)(\Omega)) = K(t) \vecse(\Omega)$ (cf.\ Lem.~\ref{lem:2}) give
\begin{equation}\label{eq:ML-K1}
\dot{{E^{\ast}}} = -e^{-\alpha(t-t_{0})}{E^{\ast}} \matse\Bigl(K^{-1}\vecse\bigl(\sum_{k} {E^{\ast}}\on{Pr}(A_{k}({E^{\ast}}))\bigr)\Bigr).
\end{equation}
By evaluation of \eqref{eq:dynamic_programming_Hessian} at $\Omega \in \se$ and application of the $\vecse-$operation to both sides of \eqref{eq:dynamic_programming_Hessian} we obtain with Lemmas \ref{lem:2} and \ref{lem3} the following dynamics of $K$:
\begin{align}\label{eq:ML-K2}
&\dot{K}\vecse(\Omega) = \Bigl( (\tilde{\Gamma}_{\vecse({E^{\ast}}^{-1}\dot{{E^{\ast}}})}^{\ast})^{\T} K + K \tilde{\Gamma}_{\vecse({E^{\ast}}^{-1}\dot{{E^{\ast}}})}^{\ast}- e^{\alpha(t-t_{0})}  KS^{-1} K \nonumber \\
& \quad + \Bigl(e^{-\alpha(t-t_{0})}   \sum_{k} \bigl( \tilde{\Gamma}_{\vecse(\on{Pr}(A_{k}({E^{\ast}})))}  + D_{k}({E^{\ast}}) \bigr)\Bigr)  \Bigr)\vecse(\Omega) .
\end{align}
Since $\Omega$ was chosen arbitrarily we can neglect $\vecse(\Omega)$ on both sides of \eqref{eq:ML-K2}. 
A change of variables $P(t):=e^{-\alpha(t-t_{0})}K(t)^{-1}$ in \eqref{eq:ML-K1} and \eqref{eq:ML-K2} gives the ODEs \eqref{eq:final-estimate-E} and \eqref{eq:final-estimate-Q} in Theorem \ref{prop:main-result}. For brevity we omit the computations here. \hfill \qed
\end{proof}
\renewcommand*{\proofname}{Proof}

\subsection{Numerical Geometric Integration}

In order to solve the differential equations \eqref{eq:final-estimate-E} and $\eqref{eq:final-estimate-Q}$ we use   \emph{Crouch-Grossman} methods \cite{hairer2006}. 
We adapt the version for right invariant Lie groups from the standard literature by permuting the order of the factors, i.e. 
\begin{equation}\label{eq:CG-left}
E_{n+1} =E_{n}  \Exp(\tfrac{h}{2}K_{1}^{E}) \Exp(\tfrac{h}{2} K_{2}^{E})\, ,
\end{equation}
with $K_{1}^{E}:= \phi(E_{n})$ and $K_{2}^{E} := \phi(E_{n}\Exp(h K_{1}^{E}))$ and step size $h$. We use the method \eqref{eq:CG-left} to integrate equation \eqref{eq:final-estimate-E}, where ${E^{\ast}}^{-1}\phi({E^{\ast}})$ is defined by the right hand side of \eqref{eq:final-estimate-E}.
 For the integration of equation \eqref{eq:final-estimate-Q} we used a standard 2-stage Runge-Kutta schemes on $\R^{6\times 6},$ since $P$ does not lie on a non-trivial manifold. 

\section{Experiments}
\subsubsection{Preprocessing.} We computed the depth map from stereo images with \cite{yamaguchi2014efficient} and the temporal optical flow between left images by the method \cite{sun2014quantitative}. Both methods are the top ranked on the Kitti Benchmark and the code is publicly available.  To remove outliers in the flow / depth map we computed for each image on 50 points $x$ the energy $E(x):= \lVert y(x) - h(x,E) \rVert,$ and removed all points $x$ with $E(x)< \lambda,$ where we selected $\lambda$ as 80\% quantile the energy of all points.

\subsubsection{Evaluation on KITTI benchmark.}
We compare our approach with Geiger et al. \cite{geiger2011stereoscan} on the challenging KITTI benchmark. 
We evaluated the first six sequences of the KITTI benchmark. Both algorithms have as initialization the identity matrix, i.e. $E_{0} = \eins$ thus it takes some frames until the approaches converge. For this reason, we omitted the first 10 frames in the evaluations.
The translational and rotational error of our approach and Geiger et al. \cite{geiger2011stereoscan} w.r.t. ground truth ego-motion in depicted in Tab. \ref{tab:eval_error}. Usually, in the rotational component our approach works better than \cite{geiger2011stereoscan}, since we model rotations explicitly on a manifold, as depicted in Fig~\ref{fig:comparison-geiger}. However, our approach less exact in the translational component because the optical flow estimation by  \cite{sun2014quantitative} fails on some frames of sequences 1 and 2 yielding a high energy and error. 
\begin{table}[tdp]
\begin{center}
\begin{tabularx}{\textwidth}{l @{\extracolsep{\fill}} rrrrrrrrrr} \hline
sequence  & 0 & 1 & 2 & 3 & 4 & 5    \\  \hline
Geiger~\cite{geiger2011stereoscan} trans. err. (m) & \textbf{0.023} & \textbf{0.050} & \textbf{0.027} & \textbf{0.017} & \textbf{0.017} & \textbf{0.017} \\
ours trans. err. (m) & 0.027 & 0.84 & 0.060 & 0.020 & 0.024 & 0.027 \\ \hline
Geiger~\cite{geiger2011stereoscan} rot. err. (deg) & 0.27 & \textbf{0.15} & \textbf{0.26} & 0.25 & \textbf{0.10} & \textbf{0.22}  \\
ours rot. err.  (deg) & \textbf{0.19} & 0.17 & 0.28 & \textbf{0.22} & \textbf{0.10} & \textbf{0.22}  \\ \hline 
\end{tabularx}
\caption{Comparison between our approach and Geiger et al. \cite{geiger2011stereoscan} on the first six sequences of the KITTI visual odometry benchmark. Our approach is usually better than \cite{geiger2011stereoscan} in the rotational part, since we model the Lie Group explicitly, but inferior in the translation since our input data from \cite{sun2014quantitative} is often not correct, (cf.~seq.1,2).}
\label{tab:eval_error} 
\end{center}
\end{table} 
\subsubsection{Synthetic Data.} 
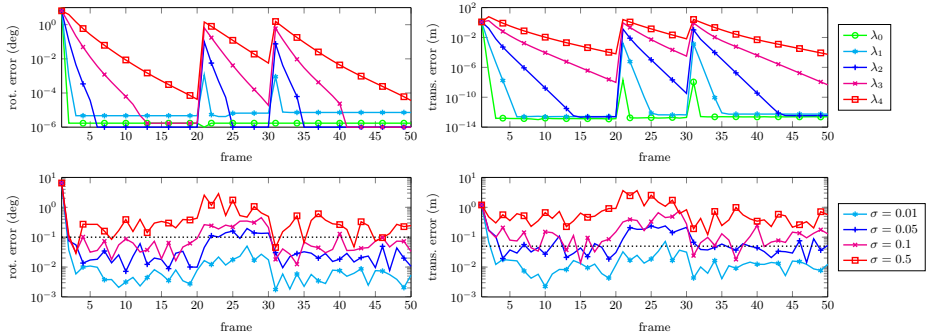
\begin{figure}[tbp] 
%
%
\begin{tikzpicture}[scale=0.52]

\begin{axis}[%
width=3.478944in,
height=1.2in,
at={(0.955556in,0.4in)},
scale only axis,
separate axis lines,
every outer x axis line/.append style={black},
every x tick label/.append style={font=\color{black}},
xmin=1,
xmax=50,
xlabel={frame},
every outer y axis line/.append style={black},
every y tick label/.append style={font=\color{black}},
ymode=log,
ymin=1e-06,
ymax=10,
ytick={ 1e-06, 0.0001,   0.01,      1},
yminorticks=true,
ylabel={rot. error (deg)}
]
\addplot [color=green,mark=o,mark options={solid}, mark repeat={3},forget plot,line width=1.0pt]
  table[row sep=crcr]{%
1	6.50754960911959\\
2	1.70754729250319e-06\\
3	1.70754729250319e-06\\
4	1.70754729250319e-06\\
5	1.70754729250319e-06\\
6	1.70754729250319e-06\\
7	1.70754729250319e-06\\
8	1.70754729250319e-06\\
9	1.70754729250319e-06\\
10	1.70754729250319e-06\\
11	1.70754729250319e-06\\
12	1.70754729250319e-06\\
13	1.70754729250319e-06\\
14	1.70754729250319e-06\\
15	1.70754729250319e-06\\
16	1.70754729250319e-06\\
17	1.70754729250319e-06\\
18	1.70754729250319e-06\\
19	1.70754729250319e-06\\
20	1.70754729250319e-06\\
21	1.e-06\\
22	1.70754729250319e-06\\
23	1.70754729250319e-06\\
24	1.70754729250319e-06\\
25	1.70754729250319e-06\\
26	1.70754729250319e-06\\
27	1.70754729250319e-06\\
28	1.70754729250319e-06\\
29	1.70754729250319e-06\\
30	1.70754729250319e-06\\
31	1.70754729250319e-06\\
32	1.70754729250319e-06\\
33	1.70754729250319e-06\\
34	1.70754729250319e-06\\
35	1.70754729250319e-06\\
36	1.70754729250319e-06\\
37	1.70754729250319e-06\\
38	1.70754729250319e-06\\
39	1.70754729250319e-06\\
40	1.70754729250319e-06\\
41	1.70754729250319e-06\\
42	1.70754729250319e-06\\
43	1.70754729250319e-06\\
44	1.70754729250319e-06\\
45	1.70754729250319e-06\\
46	1.70754729250319e-06\\
47	1.70754729250319e-06\\
48	1.70754729250319e-06\\
49	1.70754729250319e-06\\
50	1.70754729250319e-06\\
};
\addplot [color=cyan,mark=asterisk,mark options={solid}, mark repeat={3},forget plot,line width=1.0pt]
  table[row sep=crcr]{%
1	6.50754960911959\\
2	0.000939005081872024\\
3	4.67631085055434e-06\\
4	4.67631085055434e-06\\
5	4.67631085055434e-06\\
6	4.67631085055434e-06\\
7	4.67631085055434e-06\\
8	4.67631085055434e-06\\
9	4.67631085055434e-06\\
10	4.67631085055434e-06\\
11	4.67631085055434e-06\\
12	4.67631085055434e-06\\
13	4.67631085055434e-06\\
14	4.67631085055434e-06\\
15	4.67631085055434e-06\\
16	4.67631085055434e-06\\
17	4.67631085055434e-06\\
18	4.67631085055434e-06\\
19	4.67631085055434e-06\\
20	4.67631085055434e-06\\
21	0.00116164180016562\\
22	5.39973865676383e-06\\
23	4.18261957830375e-06\\
24	4.18261957830375e-06\\
25	6.50214637361826e-06\\
26	6.61330222672642e-06\\
27	6.61330222672642e-06\\
28	6.61330222672642e-06\\
29	6.61330222672642e-06\\
30	6.61330222672642e-06\\
31	0.000895451052802767\\
32	7.63638364166569e-06\\
33	7.24450961835441e-06\\
34	7.24450961835441e-06\\
35	7.24450961835441e-06\\
36	7.24450961835441e-06\\
37	7.24450961835441e-06\\
38	7.24450961835441e-06\\
39	7.24450961835441e-06\\
40	7.24450961835441e-06\\
41	7.24450961835441e-06\\
42	7.24450961835441e-06\\
43	7.24450961835441e-06\\
44	7.24450961835441e-06\\
45	7.24450961835441e-06\\
46	7.24450961835441e-06\\
47	7.24450961835441e-06\\
48	7.24450961835441e-06\\
49	7.24450961835441e-06\\
50	7.24450961835441e-06\\
};
\addplot [color=blue,mark=+,mark options={solid}, mark repeat={3},forget plot,line width=1.0pt]
  table[row sep=crcr]{%
1	6.50754960911959\\
2	0.116735377640903\\
3	0.00450255873210348\\
4	0.000341402720114882\\
5	3.65031918406285e-05\\
6	1.e-06\\
7	1.e-06\\
8	1.e-06\\
9	1.e-06\\
10	1.e-06\\
11	1.e-06\\
12	1.e-06\\
13	1.e-06\\
14	1.e-06\\
15	1.e-06\\
16	1.e-06\\
17	1.e-06\\
18	1.e-06\\
19	1.e-06\\
20	1.e-06\\
21	0.0996576327466519\\
22	0.00565808337315534\\
23	0.000499354681788403\\
24	5.60117157887658e-05\\
25	1.e-06\\
26	1.e-06\\
27	1.e-06\\
28	1.e-06\\
29	1.e-06\\
30	1.e-06\\
31	0.0755462812127567\\
32	0.00419518724594796\\
33	0.000492383926853949\\
34	6.20147701789049e-05\\
35	1.e-06\\
36	1.e-06\\
37	1.e-06\\
38	1.e-06\\
39	1.e-06\\
40	1.e-06\\
41	1.e-06\\
42	1.e-06\\
43	1.e-06\\
44	1.e-06\\
45	1.e-06\\
46	1.e-06\\
47	1.e-06\\
48	1.e-06\\
49	1.e-06\\
50	1.e-06\\
};
\addplot [color=magenta,mark=x,mark options={solid}, mark repeat={3},forget plot,line width=1.0pt]
  table[row sep=crcr]{%
1	6.50754960911959\\
2	1.22698746834983\\
3	0.214584340864831\\
4	0.0492956281236892\\
5	0.0127670518859318\\
6	0.00370209937620539\\
7	0.0012057460448812\\
8	0.000430572871657003\\
9	0.000162306838816783\\
10	6.1791038198586e-05\\
11	2.18339152038308e-05\\
12	4.82967307890293e-06\\
13	1.70754729250319e-06\\
14	1.70754729250319e-06\\
15	1.70754729250319e-06\\
16	1.70754729250319e-06\\
17	1.70754729250319e-06\\
18	1.70754729250319e-06\\
19	1.70754729250319e-06\\
20	1.70754729250319e-06\\
21	0.650673376798792\\
22	0.174159700265312\\
23	0.0478080749367776\\
24	0.013664900383118\\
25	0.00411051678392278\\
26	0.00131373655166565\\
27	0.000446547286952011\\
28	0.00015952549819876\\
29	5.79560769468377e-05\\
30	1.90909591041643e-05\\
31	0.598883688637697\\
32	0.136906114805816\\
33	0.0344702353616248\\
34	0.00945549276440521\\
35	0.00298454477969326\\
36	0.00108441151472771\\
37	0.000425635063241415\\
38	0.000170805947987499\\
39	6.73886519020876e-05\\
40	2.37220051204374e-05\\
41	1.e-06\\
42	1.e-06\\
43	1.e-06\\
44	1.e-06\\
45	1.e-06\\
46	1.e-06\\
47	1.e-06\\
48	1.e-06\\
49	1.e-06\\
50	1.e-06\\
};
\addplot [color=red,mark=square,mark options={solid}, mark repeat={3},forget plot,line width=1.0pt]
  table[row sep=crcr]{%
1	6.50754960911959\\
2	3.94509680204229\\
3	1.45694899798705\\
4	0.593383512852916\\
5	0.266136271140272\\
6	0.126350848343142\\
7	0.0621860210634532\\
8	0.0314325836428592\\
9	0.0162643353998536\\
10	0.0086140318326294\\
11	0.00467482526367423\\
12	0.00260193902680875\\
13	0.00148477594840038\\
14	0.000867006191941339\\
15	0.000516346038643453\\
16	0.000312313153895613\\
17	0.000190982122907889\\
18	0.000117362049573047\\
19	7.1574549230364e-05\\
20	4.27910125009078e-05\\
21	1.43515881866239\\
22	0.801154076926568\\
23	0.427163843760734\\
24	0.22770230112724\\
25	0.121746922053417\\
26	0.0653209374196508\\
27	0.035176933395212\\
28	0.0190195072307832\\
29	0.0103279761252223\\
30	0.00563862910174716\\
31	1.53399911006455\\
32	0.745663215781729\\
33	0.357638087627025\\
34	0.179008184370965\\
35	0.091699580683037\\
36	0.0477368899321737\\
37	0.0252824454177892\\
38	0.013688788776281\\
39	0.00761838782762971\\
40	0.00437214975206774\\
41	0.00258523421396692\\
42	0.00156788597287932\\
43	0.000969130496551597\\
44	0.000606705754389178\\
45	0.000382677310932361\\
46	0.000242258178096117\\
47	0.000153422910431564\\
48	9.68421714059033e-05\\
49	6.05276743364878e-05\\
50	3.64232285642206e-05\\
};
\end{axis}
\end{tikzpicture}%
\hspace{-0.1cm}
%
%
\begin{tikzpicture}[scale=0.52]

\begin{axis}[%
width=3.446667in,
height=1.2in,
at={(1.033333in,0.4in)},
scale only axis,
separate axis lines,
every outer x axis line/.append style={black},
every x tick label/.append style={font=\color{black}},
xmin=1,
xmax=50,
xlabel={frame},
every outer y axis line/.append style={black},
every y tick label/.append style={font=\color{black}},
ymode=log,
ymin=1e-14,
ymax=100,
ytick={1e-14, 1e-10, 1e-06,  0.01,   100},
yminorticks=true,
ylabel={trans. error (m)},
legend style={at={(1.03,0.5)},anchor=west,legend cell align=left,align=left,draw=black}
]
\addplot [color=green,mark=o,mark options={solid}, mark repeat={3},line width=1.0pt]
  table[row sep=crcr]{%
1	1.2\\
2	1.77323479069129e-07\\
3	1.48185952179579e-13\\
4	1.58446682176487e-13\\
5	1.43781987331309e-13\\
6	1.40748428625816e-13\\
7	1.38308160907782e-13\\
8	1.24750734388914e-13\\
9	1.01629281052561e-13\\
10	1.59705169123815e-13\\
11	1.38696448875984e-13\\
12	1.34001080655037e-13\\
13	1.36726449842127e-13\\
14	1.3942057469823e-13\\
15	1.34412211302483e-13\\
16	1.34921310916457e-13\\
17	1.36796847176002e-13\\
18	1.35329731692461e-13\\
19	1.23030478698153e-13\\
20	1.35297928905174e-13\\
21	2.00880168639346e-08\\
22	1.75113177385095e-13\\
23	1.67245910841189e-13\\
24	1.76354295255622e-13\\
25	1.86384200038393e-13\\
26	1.63388952802421e-13\\
27	1.71924339441437e-13\\
28	1.95187982936258e-13\\
29	1.62522966963908e-13\\
30	1.97532478119297e-13\\
31	1.09931538743351e-08\\
32	2.57710289473598e-13\\
33	2.35767718803713e-13\\
34	2.37435733586635e-13\\
35	2.48360005506057e-13\\
36	2.42197737693555e-13\\
37	2.4360693541605e-13\\
38	2.33336512577938e-13\\
39	2.35332798854857e-13\\
40	2.2869909078965e-13\\
41	2.46307542179981e-13\\
42	2.29702349863228e-13\\
43	2.40499519839503e-13\\
44	2.34890366446911e-13\\
45	2.42659853838747e-13\\
46	2.36325284325906e-13\\
47	2.41157116644315e-13\\
48	2.45256471754814e-13\\
49	2.31717511235574e-13\\
50	2.47635478835577e-13\\
};
\addlegendentry{$\lambda_0$};

\addplot [color=cyan,mark=asterisk,mark options={solid}, mark repeat={3},line width=1.0pt]
  table[row sep=crcr]{%
1	1.2\\
2	0.00144877020340833\\
3	4.49992657667179e-06\\
4	1.72939575274553e-08\\
5	6.66471663093452e-11\\
6	2.38530068305029e-13\\
7	2.55656979737888e-13\\
8	2.64431983585847e-13\\
9	2.85089530582743e-13\\
10	2.45015188301878e-13\\
11	2.55307659932684e-13\\
12	2.59988294100822e-13\\
13	2.56837563178129e-13\\
14	2.52509883652672e-13\\
15	2.59528117584662e-13\\
16	2.57438677593437e-13\\
17	2.56673409324848e-13\\
18	2.57448827973316e-13\\
19	2.73991862063046e-13\\
20	2.47943287415831e-13\\
21	0.00170392482895443\\
22	5.91510201650305e-06\\
23	2.55443980010039e-08\\
24	1.10986360303721e-10\\
25	7.65001138086911e-13\\
26	4.4168609412156e-13\\
27	4.4180534131167e-13\\
28	4.23437197581013e-13\\
29	4.63606437517226e-13\\
30	4.31074111364252e-13\\
31	0.00133005051029216\\
32	5.81112935870513e-06\\
33	2.24595835806184e-08\\
34	8.66357864610836e-11\\
35	8.01751091912009e-13\\
36	5.46274967907888e-13\\
37	5.45732955071867e-13\\
38	5.59597402299784e-13\\
39	5.56660321387929e-13\\
40	5.63632742787833e-13\\
41	5.53097780787306e-13\\
42	5.58833233044291e-13\\
43	5.55405667293914e-13\\
44	5.61039521773529e-13\\
45	5.48982549230793e-13\\
46	5.62655838660767e-13\\
47	5.57393282646942e-13\\
48	5.53819722674601e-13\\
49	5.64678717725193e-13\\
50	5.48223695090601e-13\\
};
\addlegendentry{$\lambda_1$};

\addplot [color=blue,mark=+,mark options={solid}, mark repeat={3},line width=1.0pt]
  table[row sep=crcr]{%
1	1.2\\
2	0.168163949525161\\
3	0.00657932816169094\\
4	0.000516948101867929\\
5	5.75092925464252e-05\\
6	6.98427614274919e-06\\
7	8.64925536853745e-07\\
8	1.07604288723883e-07\\
9	1.3401652934941e-08\\
10	1.66953675793805e-09\\
11	2.07909228962366e-10\\
12	2.57934118467773e-11\\
13	3.1010978766165e-12\\
14	3.37995621944024e-13\\
15	2.1124049807545e-13\\
16	2.39939518371919e-13\\
17	2.38964888907694e-13\\
18	2.41059055065519e-13\\
19	2.44657614779923e-13\\
20	2.51250327075208e-13\\
21	0.136044111474764\\
22	0.00797843693608831\\
23	0.000721465982814457\\
24	8.30457676849993e-05\\
25	1.03046648085991e-05\\
26	1.30346691951264e-06\\
27	1.65372676657022e-07\\
28	2.09298725580615e-08\\
29	2.63622655034185e-09\\
30	3.30269260874336e-10\\
31	0.103146436168814\\
32	0.00570638099263404\\
33	0.000731412129043791\\
34	9.5340074693222e-05\\
35	1.20069744669798e-05\\
36	1.49415630337595e-06\\
37	1.85368876842602e-07\\
38	2.2994481933859e-08\\
39	2.85449031688158e-09\\
40	3.54690269987811e-10\\
41	4.41179097365358e-11\\
42	5.49967995729303e-12\\
43	7.57377444631379e-13\\
44	3.50381224321943e-13\\
45	3.45925631889854e-13\\
46	3.45563682593476e-13\\
47	3.43915341937582e-13\\
48	3.47246722057801e-13\\
49	3.52629795790057e-13\\
50	3.51394999622199e-13\\
};
\addlegendentry{$\lambda_2$};

\addplot [color=magenta,mark=x,mark options={solid}, mark repeat={3},line width=1.0pt]
  table[row sep=crcr]{%
1	1.2\\
2	1.81788512943983\\
3	0.291726696665035\\
4	0.0645275167856574\\
5	0.0170815955289867\\
6	0.00514986908715237\\
7	0.00172538652398052\\
8	0.00062370080164347\\
9	0.000236283576407978\\
10	9.19117677884563e-05\\
11	3.62742795734363e-05\\
12	1.44306379487702e-05\\
13	5.76636259315326e-06\\
14	2.30999648499371e-06\\
15	9.2671897806913e-07\\
16	3.7209122776488e-07\\
17	1.49474045757816e-07\\
18	6.00635912614039e-08\\
19	2.41399621403028e-08\\
20	9.70315539957497e-09\\
21	0.92540786457851\\
22	0.220304567226303\\
23	0.0593241443488504\\
24	0.0172185932673064\\
25	0.0053213469273468\\
26	0.00174917994124604\\
27	0.000608468292509379\\
28	0.000221460926244894\\
29	8.32052310816369e-05\\
30	3.19116971884445e-05\\
31	0.875425478319718\\
32	0.173901450310149\\
33	0.040216984822441\\
34	0.0109330093433004\\
35	0.00367994803867607\\
36	0.00143951009643545\\
37	0.000591429690114013\\
38	0.000243923877406747\\
39	9.99284640526265e-05\\
40	4.06591786124519e-05\\
41	1.64623295244481e-05\\
42	6.64454391664552e-06\\
43	2.67698390499619e-06\\
44	1.07748203483515e-06\\
45	4.33510626707956e-07\\
46	1.74408227616025e-07\\
47	7.01788746446094e-08\\
48	2.82472822218946e-08\\
49	1.1374037558406e-08\\
50	4.58187158849291e-09\\
};
\addlegendentry{$\lambda_3$};

\addplot [color=red,mark=square,mark options={solid}, mark repeat={3},line width=1.0pt]
  table[row sep=crcr]{%
1	1.2\\
2	6.04011959134077\\
3	2.15273718157215\\
4	0.82088479435853\\
5	0.345565777390055\\
6	0.157276150478263\\
7	0.0759871926781798\\
8	0.0384328658174683\\
9	0.0201595601765932\\
10	0.0109038285018574\\
11	0.00605950299470289\\
12	0.00345026231819613\\
13	0.00200734307396622\\
14	0.00118971086474355\\
15	0.000716074730857491\\
16	0.000436391186031581\\
17	0.000268560809647341\\
18	0.000166530856973554\\
19	0.000103862830087928\\
20	6.50636180934786e-05\\
21	2.53469299532189\\
22	1.16621495162685\\
23	0.553192939331532\\
24	0.275257280593429\\
25	0.141863333397741\\
26	0.0748586159198897\\
27	0.0401304646379831\\
28	0.0217545671765958\\
29	0.0118945638500789\\
30	0.00655482025049716\\
31	2.55699436147666\\
32	1.10452722794743\\
33	0.473478404459357\\
34	0.215229919000885\\
35	0.102829060783544\\
36	0.0513044627224726\\
37	0.0267394771806684\\
38	0.0146099619045578\\
39	0.00837723085752731\\
40	0.00501327757695068\\
41	0.00309779567139643\\
42	0.00195375334898953\\
43	0.00124618240828059\\
44	0.000799028255736974\\
45	0.000513203345034884\\
46	0.000329587036264583\\
47	0.000211466133437874\\
48	0.000135511116899362\\
49	8.67296539392477e-05\\
50	5.54461561043295e-05\\
};
\addlegendentry{$\lambda_4$};

\end{axis}
\end{tikzpicture}%
 \\
%
%
\begin{tikzpicture}[scale=0.52]

\begin{axis}[%
width=3.478944in,
height=1.2in,
at={(0.955556in,0.4in)},
scale only axis,
separate axis lines,
every outer x axis line/.append style={black},
every x tick label/.append style={font=\color{black}},
xmin=1,
xmax=50,
xlabel={frame},
every outer y axis line/.append style={black},
every y tick label/.append style={font=\color{black}},
ymode=log,
ymin=1e-03,
ymax=10,
ytick={ 1e-3,  1e-2,   1e-1, 1, 10},
yminorticks=true,
ylabel={rot. error (deg)}
]
\addplot [color=cyan,solid,mark=asterisk,mark options={solid},mark repeat={3},forget plot,line width=1.0pt]
  table[row sep=crcr]{%
1	6.50754960911959\\
2	0.0370251838784492\\
3	0.00589173722476192\\
4	0.0100286215128469\\
5	0.0108437053891719\\
6	0.00984671638456139\\
7	0.00377224946503588\\
8	0.00383014295815885\\
9	0.0020868858238938\\
10	0.00316570334174755\\
11	0.00451067387778668\\
12	0.00233302112230087\\
13	0.00445309837366231\\
14	0.00819008456493369\\
15	0.00575525566499569\\
16	0.00352709569258308\\
17	0.00799624552668476\\
18	0.00347808522045975\\
19	0.00274061606445858\\
20	0.00782974152750454\\
21	0.016498735668377\\
22	0.0121795995942845\\
23	0.0216887444405236\\
24	0.0116788029007747\\
25	0.0176656883322904\\
26	0.0230412994419543\\
27	0.0502680098601407\\
28	0.0169413010787455\\
29	0.0268703390766808\\
30	0.0174921654969797\\
31	0.00176474124628178\\
32	0.00586575501912273\\
33	0.00184225717663988\\
34	0.00587462556690716\\
35	0.00226097123019296\\
36	0.00566287303943643\\
37	0.00346046034026948\\
38	0.0028144223174438\\
39	0.00789064662270596\\
40	0.00759394526753877\\
41	0.00492819773676244\\
42	0.00620031031413528\\
43	0.00874217114050574\\
44	0.00247434797549252\\
45	0.00761260644968746\\
46	0.00735816831775831\\
47	0.00559586904066162\\
48	0.00388421765863078\\
49	0.00204823340560525\\
50	0.00539522002170967\\
};
\addplot [color=blue,solid,mark=+,mark options={solid},mark repeat={3},forget plot,line width=1.0pt]
  table[row sep=crcr]{%
1	6.50754960911959\\
2	0.0754857703618961\\
3	0.05546093668958\\
4	0.0136119689005089\\
5	0.0172258585989981\\
6	0.0182514292162717\\
7	0.0344008110217811\\
8	0.0138776674371293\\
9	0.0255957191980776\\
10	0.00717956367416877\\
11	0.0229944135759313\\
12	0.0606239799033748\\
13	0.0122961579252318\\
14	0.0200353147980905\\
15	0.0179486628969276\\
16	0.00886951463780675\\
17	0.0291307303746133\\
18	0.0336956300951751\\
19	0.0100595791505567\\
20	0.00992024448046902\\
21	0.045338292735864\\
22	0.112855224805392\\
23	0.101938843149399\\
24	0.140839418070186\\
25	0.159267009714735\\
26	0.0911928656628804\\
27	0.195543188744992\\
28	0.138821016857337\\
29	0.134113734978335\\
30	0.137075663338712\\
31	0.0305266684507302\\
32	0.00953979456803569\\
33	0.022255292769784\\
34	0.0284189301586329\\
35	0.0325014418755793\\
36	0.0292111107913681\\
37	0.0198716021296267\\
38	0.0285211240781952\\
39	0.0153643229006501\\
40	0.0184050794309422\\
41	0.0287071617235372\\
42	0.0147907201819551\\
43	0.0207572158948542\\
44	0.0233080128849433\\
45	0.00962610564876895\\
46	0.0494574112198113\\
47	0.00767621955725361\\
48	0.0436896842345416\\
49	0.0198644943509393\\
50	0.0105800786481897\\
};
\addplot [color=magenta,solid,mark=x,mark options={solid}, mark repeat={3}, forget plot,line width=1.0pt]
  table[row sep=crcr]{%
1	6.50754960911959\\
2	0.0885967682885412\\
3	0.0312266895287464\\
4	0.10231637933185\\
5	0.0315256532230513\\
6	0.0361507207686293\\
7	0.0736625882810716\\
8	0.0180278238964643\\
9	0.0551326454680579\\
10	0.0869697628726404\\
11	0.0642346791213029\\
12	0.031802017689728\\
13	0.0373488111091715\\
14	0.0730528266381397\\
15	0.018974359798615\\
16	0.0645247464545221\\
17	0.0163011089355162\\
18	0.0468346224968836\\
19	0.0621928089866154\\
20	0.0909195841190389\\
21	0.262898147753191\\
22	0.243323788606644\\
23	0.188425978133381\\
24	0.253266169117915\\
25	0.221536620862377\\
26	0.3435578493184\\
27	0.346914400815566\\
28	0.307954857456214\\
29	0.451281441806276\\
30	0.211309820423905\\
31	0.0180767779693975\\
32	0.0405179864167133\\
33	0.0271140533005113\\
34	0.0124815670033702\\
35	0.069260133665313\\
36	0.0501749809303349\\
37	0.0447961960997471\\
38	0.0301686465006484\\
39	0.0389866960790554\\
40	0.128883848934549\\
41	0.0357647371069428\\
42	0.0365652350997377\\
43	0.0426238451519928\\
44	0.0447412816240542\\
45	0.0735586108749001\\
46	0.0977732913701167\\
47	0.0472247410246389\\
48	0.0735170259201252\\
49	0.072773272951309\\
50	0.0333333304960017\\
};
\addplot [color=red,solid,mark=square,mark options={solid}, mark repeat={3},forget plot,line width=1.0pt]
  table[row sep=crcr]{%
1	6.50754960911959\\
2	0.116948475106023\\
3	0.0287923948486935\\
4	0.271733297098814\\
5	0.269096704585752\\
6	0.269281864031831\\
7	0.159635846609703\\
8	0.0846751299817035\\
9	0.166419532996169\\
10	0.388878058491996\\
11	0.155429098185811\\
12	0.393564528287142\\
13	0.142135279194413\\
14	0.296308148384351\\
15	0.3366290516042\\
16	0.296748551503562\\
17	0.205436536041431\\
18	0.390899245165891\\
19	0.3769231965623\\
20	0.432379790563057\\
21	2.54981432176983\\
22	1.17035801380185\\
23	2.85837984027745\\
24	0.534229094554339\\
25	1.75131749346836\\
26	0.528696502397351\\
27	0.468824941271289\\
28	1.06337440652121\\
29	0.649991968706905\\
30	0.495901518242898\\
31	0.0453156141067055\\
32	0.153172775839503\\
33	0.239620891541756\\
34	0.519784882240313\\
35	0.0664703565638408\\
36	0.36280012156607\\
37	0.608036274353685\\
38	0.372970455247856\\
39	0.280826920310875\\
40	0.264321992296242\\
41	0.181072776978737\\
42	0.0474788649518271\\
43	0.329142650792365\\
44	0.245686138737216\\
45	0.0741949096215426\\
46	0.0965794476192978\\
47	0.132121470366061\\
48	0.269294015787893\\
49	0.225629967915626\\
50	0.24516813133425\\
};
\addplot [color=black,dotted,forget plot,line width=1.pt]
  table[row sep=crcr]{%
1	0.1\\
2	0.1\\
3	0.1\\
4	0.1\\
5	0.1\\
6	0.1\\
7	0.1\\
8	0.1\\
9	0.1\\
10	0.1\\
11	0.1\\
12	0.1\\
13	0.1\\
14	0.1\\
15	0.1\\
16	0.1\\
17	0.1\\
18	0.1\\
19	0.1\\
20	0.1\\
21	0.1\\
22	0.1\\
23	0.1\\
24	0.1\\
25	0.1\\
26	0.1\\
27	0.1\\
28	0.1\\
29	0.1\\
30	0.1\\
31	0.1\\
32	0.1\\
33	0.1\\
34	0.1\\
35	0.1\\
36	0.1\\
37	0.1\\
38	0.1\\
39	0.1\\
40	0.1\\
41	0.1\\
42	0.1\\
43	0.1\\
44	0.1\\
45	0.1\\
46	0.1\\
47	0.1\\
48	0.1\\
49	0.1\\
50	0.1\\
};
\end{axis}
\end{tikzpicture}%
\hspace{-0.1cm}
%
%
\begin{tikzpicture}[scale=0.52]

\begin{axis}[%
width=3.446667in,
height=1.2in,
at={(1.033333in,0.4in)},
scale only axis,
separate axis lines,
every outer x axis line/.append style={black},
every x tick label/.append style={font=\color{black}},
xmin=1,
xmax=50,
xlabel={frame},
every outer y axis line/.append style={black},
every y tick label/.append style={font=\color{black}},
ymode=log,
ymin=1e-3,
ymax=10,
ytick={1e-3, 1e-2, 1e-01,  1,   10},
yminorticks=true,
ylabel={trans. error (m)},
legend style={at={(1.03,0.5)},anchor=west,legend cell align=left,align=left,draw=black}
]

\addplot [color=cyan,solid,mark=asterisk,mark options={solid},mark repeat={3},line width=1.0pt]
  table[row sep=crcr]{%
1	1.2\\
2	0.0690882786543981\\
3	0.0122805995357967\\
4	0.0179378127819187\\
5	0.0169286929769643\\
6	0.0161297405743997\\
7	0.00744576053020204\\
8	0.00424084428109044\\
9	0.00575269718997841\\
10	0.00223919361213584\\
11	0.00585407072083983\\
12	0.0083243354480555\\
13	0.00821225654250034\\
14	0.0165310714943838\\
15	0.0132198944792479\\
16	0.00950529932365631\\
17	0.012558215198977\\
18	0.00606509246539031\\
19	0.00429561122200703\\
20	0.0117200079974794\\
21	0.0325390871202778\\
22	0.0168601125780953\\
23	0.0247763399597177\\
24	0.0135373394187872\\
25	0.0186220708470677\\
26	0.0335430525460111\\
27	0.0719309325562703\\
28	0.0166564207737687\\
29	0.036714864917179\\
30	0.0315876122225964\\
31	0.00430855348705766\\
32	0.0100764786184378\\
33	0.00621099349366018\\
34	0.0165023403280652\\
35	0.00409693697578698\\
36	0.0110793797832357\\
37	0.00911160225599947\\
38	0.0122455606633535\\
39	0.0146229845451981\\
40	0.0124221137996283\\
41	0.0123133858581265\\
42	0.0142698754846713\\
43	0.0133708198785033\\
44	0.00519302790343959\\
45	0.0158732814189576\\
46	0.0152346359108969\\
47	0.0123544663057103\\
48	0.00871046513475687\\
49	0.00769236555478883\\
50	0.0124092334101271\\
};
\addlegendentry{$\sigma=0.01$};

\addplot [color=blue,solid,mark=+,mark options={solid},mark repeat={3},line width=1.0pt]
  table[row sep=crcr]{%
1	1.2\\
2	0.172436670183728\\
3	0.106507055405071\\
4	0.0193155323262614\\
5	0.0373973455963369\\
6	0.0294167596272366\\
7	0.0567293496777467\\
8	0.0319846422105164\\
9	0.0451680037334791\\
10	0.0270539622054182\\
11	0.0537194069450236\\
12	0.109023238120186\\
13	0.021062434004111\\
14	0.0286246377828296\\
15	0.0436133797904131\\
16	0.0222977744330694\\
17	0.0610079261109765\\
18	0.0461570487373744\\
19	0.0186115390100738\\
20	0.0306149268920846\\
21	0.0876878990407253\\
22	0.209452291843042\\
23	0.18341292330426\\
24	0.184394293476628\\
25	0.236254420309675\\
26	0.200558336472054\\
27	0.238007711680474\\
28	0.141742156929724\\
29	0.175166120369936\\
30	0.212970752300814\\
31	0.0676083189919541\\
32	0.0387475591005216\\
33	0.0374744692615725\\
34	0.0455400434234328\\
35	0.0588981137856288\\
36	0.0460548530796427\\
37	0.0370209229810136\\
38	0.0606445657577664\\
39	0.0495192628538178\\
40	0.0442645998619698\\
41	0.0329927006618611\\
42	0.020771877315175\\
43	0.0337368594461129\\
44	0.0327060943158918\\
45	0.029488187939899\\
46	0.0758504674053407\\
47	0.0142546625072258\\
48	0.0883731039052643\\
49	0.0366938888962516\\
50	0.0530362832010125\\
};
\addlegendentry{$\sigma=0.05$};

\addplot [color=magenta,solid,mark=x,mark options={solid},mark repeat={3},line width=1.0pt]
  table[row sep=crcr]{%
1	1.2\\
2	0.16781814850807\\
3	0.0977020068496079\\
4	0.208325071813729\\
5	0.0796370768823628\\
6	0.074811550809373\\
7	0.144104052349886\\
8	0.0444115752218255\\
9	0.0342101899479657\\
10	0.151242880406785\\
11	0.141876343339874\\
12	0.0983489811056148\\
13	0.0374164385215752\\
14	0.113827266663538\\
15	0.0138227840584138\\
16	0.0952175729178047\\
17	0.0322958645845113\\
18	0.0711246775823012\\
19	0.0855818277813361\\
20	0.150975933779913\\
21	0.321940658699789\\
22	0.391940749893418\\
23	0.362078119578076\\
24	0.166340453181402\\
25	0.332730902320661\\
26	0.660886937975424\\
27	0.453534361985545\\
28	0.487241710667271\\
29	0.777234944759366\\
30	0.278795332094298\\
31	0.024153099318373\\
32	0.0653461815613446\\
33	0.0987054340103151\\
34	0.0533536405909864\\
35	0.130113461672676\\
36	0.0798836587137725\\
37	0.0180704235762425\\
38	0.0551401065882279\\
39	0.106173334321361\\
40	0.246862117682527\\
41	0.048552662232413\\
42	0.0686715815173906\\
43	0.0774195650230845\\
44	0.0667229673895623\\
45	0.147414840944684\\
46	0.13632902747926\\
47	0.0952145489919175\\
48	0.135014554184494\\
49	0.178065660258341\\
50	0.129409490237872\\
};
\addlegendentry{$\sigma=0.1$};

\addplot [color=red,solid,mark=square,mark options={solid},mark repeat={3},line width=1.0pt]
  table[row sep=crcr]{%
1	1.2\\
2	0.47917637993039\\
3	0.300546993472625\\
4	0.563378099269353\\
5	0.359150184853485\\
6	0.611145138038804\\
7	0.522405991844895\\
8	0.229063735427917\\
9	0.353339227329089\\
10	0.666515291446746\\
11	0.426597018062263\\
12	0.725840147434513\\
13	0.225436099153257\\
14	0.349666020214876\\
15	0.640761930664612\\
16	0.462567495471052\\
17	0.504700835968022\\
18	0.745432326625011\\
19	0.909087575631874\\
20	0.914980253837503\\
21	3.62907427638266\\
22	2.44610510587268\\
23	3.6303773968131\\
24	1.12246743154469\\
25	2.54000780581416\\
26	1.4951115416662\\
27	1.07965861982509\\
28	1.73205595902142\\
29	0.616719182075088\\
30	0.829449776312911\\
31	0.192159140383712\\
32	0.424307465034648\\
33	0.122098119748636\\
34	0.727030338884643\\
35	0.245326674125833\\
36	0.377189294192436\\
37	0.988193880937017\\
38	0.559283019981073\\
39	0.408553010042859\\
40	0.34353179899335\\
41	0.380087454885221\\
42	0.22742863982398\\
43	0.586379059735721\\
44	0.517700906192603\\
45	0.21508127255309\\
46	0.261711911350039\\
47	0.330104450548265\\
48	0.203802177158284\\
49	0.713962846994526\\
50	0.565101371264091\\
};
\addlegendentry{$\sigma=0.5$};

\addplot [color=black,dotted,forget plot,line width=1.pt]
  table[row sep=crcr]{%
1	0.05\\
2	0.05\\
3	0.05\\
4	0.05\\
5	0.05\\
6	0.05\\
7	0.05\\
8	0.05\\
9	0.05\\
10	0.05\\
11	0.05\\
12	0.05\\
13	0.05\\
14	0.05\\
15	0.05\\
16	0.05\\
17	0.05\\
18	0.05\\
19	0.05\\
20	0.05\\
21	0.05\\
22	0.05\\
23	0.05\\
24	0.05\\
25	0.05\\
26	0.05\\
27	0.05\\
28	0.05\\
29	0.05\\
30	0.05\\
31	0.05\\
32	0.05\\
33	0.05\\
34	0.05\\
35	0.05\\
36	0.05\\
37	0.05\\
38	0.05\\
39	0.05\\
40	0.05\\
41	0.05\\
42	0.05\\
43	0.05\\
44	0.05\\
45	0.05\\
46	0.05\\
47	0.05\\
48	0.05\\
49	0.05\\
50	0.05\\
};
\end{axis}
\end{tikzpicture}%
 \\
\vspace{-2.5em}
\caption{Evaluation of our method on synthetic data. In the top row we observe the linear convergence in the logarithmic scale, which means exponential convergence behavior. In frames 21 and 31 there is an immediate change of direction which causes errors since the motion-constancy assumption is violated. However, our method adapts to the change in the subsequent steps. In the bottom row we consider multiplicative Gaussian noise with mean 1 and standard deviation $\sigma$ to the input data (optical flow). For little noise rates ($\sigma\le 0.01$) we obtain the accuracy required for practical applications (dotted lines): ($0.1$ degree, $0.05$ meters).}
\label{fig:synthetic-data}
\end{figure}

We evaluated our method on synthetic data with known depth maps and optical flow. Fig.~\ref{fig:synthetic-data} (top) shows the convergence from a wrong initialization for different weights of the penalty term of $\delta,$ i.e. $S_{i} = \lambda_{i}\diag(s_{1},s_{1},s_{1},s_{2},s_{2},s_{2})$ with $s_{1} = 10^{-3},\,s_{2}=10^{-6}$ and $\lambda_{i}=10^{i}$ for $i=0,\dots,4$, $\alpha=0$ and $Q=0.1\eins.$  In frames 21 and 51 constant motion assumption is violated, leading to a high error. However, for small weights of the penalty term of~$\delta$ the filter converges almost immediately. Fig.~\ref{fig:synthetic-data} (bottom) shows the performance of the filter on data distorted by multiplicative Gaussian noise. For high noise rates ($\sigma>0.1$) the filter fails while for small noise rates ($\sigma<0.05$) filter results have an accuracy comparable to state-of-the-art filters on real data. On the other hand this means that the input data is allowed to be wrong up to 1\% in order to reach state-of-the-art results. For our evaluations on Tab.~\ref{tab:eval_error} this was not always the case.

\section{Conclusion} We presented a second order Minimum Energy Filter with non-linear observation equations for the ego-motion estimation problem and derived explicit differential equations for the optimal state. 
Our experiments showed that our approach is comparable with the state-of-the-art approaches \cite{geiger2011stereoscan}: In the translational component our method is inferior but it is  superior in the rotational component. The experiments also confirm the exponential convergence rate and robustness against multiplicative noise.

In future work, we will generalize our model by allowing an acceleration of the camera. We expect this generalization to
reduce the error in the translational part. Moreover, we plan to design a filter for monocular depth map estimation.

\vspace{0.5em}
\noindent
\textbf{Acknowledgments.}This work was supported by the DFG, grant GRK 1653. The final publication is available at link.springer.com 
\appendix
\vspace{-1em}
\section{Proofs of Lemmas~\ref{lemma:riemannian-gradient}--\ref{lem3}}
\vspace{-0.5em}
\begin{proof}[of Lem.~\ref{lemma:riemannian-gradient}]
We begin with the directional derivative of $h_{k}$ which is 
\begin{equation}
\D h_{k}(E)[E\Omega] = -\kappa_{k}^{-1} \hat{I} \Omega E^{-1} g_{k} + \kappa_{k}^{-2} (e_{3}^{\T}\hat{I} \Omega E^{-1} g_{k}) \hat{I}E^{-1} g_{k}, 
\end{equation}
where $\kappa_{k} = \kappa_{k}(E):= e_{3}^{\T}\hat{I}E^{-1} g_{k}.$  Then the following holds:
\begin{align}
\D_{1} \mathcal{H}^{-}(E, \mu, t)[E\Omega] = -e^{-\alpha(t-t_{0})} \sum_{k=1}^{n} \tr \bigl( \D h_{k}(E)[E\Omega] (y_{k} - h_{k}(E))^{\T} Q \bigr) \nonumber \\
 = e^{-\alpha(t-t_{0})} \sum_{k=1}^{n}  \big \la \bigl(\kappa_{k}^{-1}\hat{I}  -  \kappa_{k}^{-2}  \hat{I}   E^{-1}  e_{3} g_{k}^\T   \hat{I} \bigr)^{\T}Q(y_{k} - h_{k}(E))   g_{k}^{\T}  E^{-\T}  ,\Omega  \big \ra \, .\label{eq:proof_lemma1:DH}
\end{align}
We obtain the Riemannian gradient on~$\SE$ by projecting (cf.~\cite[Sec.~3.6.1]{Absil2008}) the left hand side of the Riemannian metric in \eqref{eq:proof_lemma1:DH} onto $T_{E}\SE$, which is
\begin{align*}
\grad_{1} \mathcal H^{-}(E,\mu,t) = & e^{-\alpha(t-t_{0})} \sum_{k} E\on{Pr}\bigl( A_{k}(E)  \bigr).
\end{align*}
with $A_{k}(E):=  
\bigl(\kappa_{k}^{-1}\hat{I}  -  \kappa_{k}^{-2}  \hat{I}   E^{-1}  e_{3} g_{k}^\T   \hat{I} \bigr)^{\T}
Q(y_{k} - h_{k}(E))   g_{k}^{\T}  E^{-\T}$ 
and $\on{Pr}$ is this mentioned projection.\hfill \qed
\end{proof}

\begin{proof}[of Lem.~\ref{lem:2}] The existence of a matrix~$K$ such that~$K\vecse(\Omega) = Z({E^{\ast}},t)[\Omega]$ follows from considering the basis of~$\se$ (cf.~Sec.~\ref{sec1}) and also its inverse representation.   The Hessian of the Hamiltonian w.r.t.\ the second argument evaluated at zero is $\vecse (\Hess_{2} \mathcal H^{-} ({E^{\ast}},0,t)[\Omega]) = - e^{\alpha(t-t_{0})} S^{-1}\vecse(\Omega)$, thus
$$\vecse( Z({E^{\ast}},t) (\Hess_{2}  \mathcal H^{-}({E^{\ast}},0,t)[ Z({E^{\ast}},t) (\Omega)]) ) =  - e^{\alpha(t-t_{0})} K(t)S^{-1}K(t) \vecse(\Omega).$$
 For left-invariant Lie groups the vectorized representation of affine connection holds (cf.~\cite{Absil2008}), i.e.
\begin{equation}\label{eq:flatten-levi-civita}
\vecse(\omega_{\Omega}\Delta) = \tilde{\Gamma}_{\vecse(\Delta)} \vecse(\Omega) + \vecse(\D  \Delta [\Omega]) \,,
\end{equation}
where matrix~$\tilde \Gamma_{\vecse(\Delta)}$ has the components $(\tilde \Gamma_{\vecse(\Delta)})_{ij} = \sum_{k=1}^{6}\hat\Gamma_{jk}^{i} \Delta^{k}$. Note that we yield the definition in~\cite{Absil2008} by exchanging~$i$ and~$j$. In the case of a constant function~$\Delta$
we have~$\D \Delta[\Omega]=0.$ For the dual expression it holds for all~$\Omega, \Delta \in \se$:
$\vecse(\omega^{\ast}_{\Omega}\Delta) = \tilde{\Gamma}^{\ast}_{\vecse(\Delta)} \vecse(\Omega)$
with~$(\tilde \Gamma_{\vecse(\Delta)}^{\ast})_{ik} = \sum_{j=1}^{6}\hat\Gamma_{jk}^{i} \Delta^{j}.$ \hfill \qed
\end{proof}
\begin{proof}[of Lem.~\ref{lem3}] For any matrices $A,B \in \R^{4\times4}$ and $\Omega \in \se$ 
let $\kronse, \kronse^{\T}$ denote the operators which  extract the vector form of $\Omega$ 
and are defined as $\vecse(A \Omega B)=: (A \kronse B) \vecse(\Omega)$ and $\vecse(A \Omega^{\T} B) =: (A \kronse^{\T} B)\vecse(\Omega)$, respectively, for which explicit expressions exist.
By definition of the Hessian, Lemma~\ref{lemma:riemannian-gradient} and with Eq.~\eqref{eq:flatten-levi-civita} we obtain
\begin{align*}
\vecse(& {E^{\ast}}^{-1} \Hess_{1} \mathcal H^{-}({E^{\ast}},0,t) [{E^{\ast}}\Omega]) = \vecse({E^{\ast}}^{-1} \nabla_{{E^{\ast}}\Omega} \grad_{1}\mathcal H^{-}({E^{\ast}},0,t))  \\
 = & e^{-\alpha (t-t_{0})} \sum_{k} \Bigl( \tilde{\Gamma}_{\vecse(\on{Pr}(A_{k}({E^{\ast}})))}\vecse(\Omega) +  \vecse(\mathbf D_{{E^{\ast}}} \on{Pr}(A_{k}({E^{\ast}}))[\Omega]) \Bigr) \,.
\end{align*}
It can be shown that we can omit the projection $\on{Pr}$, i.e.~$\vecse( \D \on{Pr}(A_{k}(E)) [\Omega])  = \vecse(\D (A_{k}(E))[\Omega])$ for all $\Omega \in \se.$ After computing the directional derivative of $A_{k}(E)$ in direction $\Omega$ we apply the $\vecse-$ operation and extract with the $\kronse-, \kronse^{\T}-$ operations the direction $\Omega.$ This gives
\begin{align}
\vecse(&\D A_{k}(E) [\Omega]) = 
 \Bigl(\kappa_{k}^{-2}\hat{I}^{\T}  Q (y_{k} - h_{k}(E))e_{3}^{\T} \hat{I}E^{-1} \kronse E^{-1} g_{k}g_{k}^{\T} E^{-\T}  \nonumber \\
&  - 2 \kappa_{k}^{-3} \hat{I}^{\T} e_{3}e_{3}^{\T} \hat{I}E^{-1} \kronse E^{-1} g_{k}g_{k}^{\T}E^{-\T}\hat{I}^{\T}Q (y_{k} - h_{k}(E)) g_{k}^{\T} E^{-\T}\nonumber \\
&  + \kappa_{k}^{-2} \hat{I}^{\T}e_{3}g_{k}^{\T}E^{-\T} \kronse^{\T}E^{-\T}\hat{I}^{\T}Q(y_{k} - h_{k}(E)) g_{k}^{\T} E^{-\T} \nonumber \\
& -  \bigl(\kappa_{k}^{-1}\hat{I}^{\T} -  \kappa_{k}^{-2}  \hat{I}^{\T}e_{3}g_{k}^{\T}     E^{-\T}\hat{I}^{\T} \bigr)  Q(y_{k}-h_{k}(E)) g_{k}^{\T}  E^{-\T}\kronse^{\T} E^{-\T} \nonumber \\
& -  \bigl(\kappa_{k}^{-3}\hat{I}^{\T} -  \kappa_{k}^{-4}  \hat{I}^{\T}e_{3}g_{k}^{\T}     E^{-\T}\hat{I}^{\T} \bigr)  \hat{I}E^{-1} g_{k}e_{3}^{\T}Q\hat{I}E^{-1} \kronse E^{-1} g_{k}g_{k}^{\T}E^{-\T}\nonumber \\
& +\bigl(\kappa_{k}^{-2}\hat{I}^{\T} -  \kappa_{k}^{-3}  \hat{I}^{\T}e_{3}g_{k}^{\T}     E^{-\T}\hat{I}^{\T} \bigr) Q \hat{I} E^{-1}\kronse E^{-1} g_{k}g_{k}^{\T} E^{-\T}  \Bigr)\vecse(\Omega) \nonumber \\
=:& D_{k}(E) \vecse(\Omega)  \,. \label{eq:def_Dk} 
\end{align}
\end{proof}

\bibliographystyle{plain}
\bibliography{Berger2015}

\end{document}